\documentstyle{article}
\newtheorem {th}{Theorem}[section]
\newtheorem {lem}[th]{Lemma}
\newtheorem {pr}[th]{Proposition}

\newtheorem{defn}[th]{Definition}

\newtheorem{conj}{Conjecture}
\newtheorem{qq}[conj]{Question}

\def\Cox{\hfill \Box}
\def\T{{\bf T}}
\def\had{\odot}
\def\stoch{\succeq}
\def\predeq{\preceq}

\def\diseq{\, {\stackrel {{\cal D}} {=}}}

\def\ee{\epsilon}
\def\E{{\bf{E}}}
\def\P{{\bf{P}}}
\def\N{\hbox{I\kern-.2em\hbox{N}}}
\def\R{\hbox{I\kern-.2em\hbox{R}}}

\def\B{{\cal{B}}}

\def\S{{\cal{S}}}

\def\|{\, | \, }
\def\one{{\bf 1}}
\def\0{\hat{0}}
\def\1{\hat{1}}

\def\neg{\mbox{S-MRR}_2}
\def\Cov{{\rm Cov}}
\def\ul{\underline}
\def\ngeq{\setminus \kern-1.0em \geq}

\begin{document}

\begin{titlepage}
\begin{center}
{\large \bf Towards a theory of negative dependence} \\
\end{center}
\vspace{5ex}
\begin{flushright}
Robin Pemantle \footnote{Research supported in part by 
National Science Foundation grant \# DMS 9300191, by a Sloan Foundation
Fellowship, and by a Presidential Faculty Fellowship}$^,$\footnote{Department 
of Mathematics, University of Wisconsin-Madison, Van Vleck Hall, 480 Lincoln
Drive, Madison, WI 53706.  Now at Department of Mathematics, Ohio State
University, 231 W. 18th Avenue, Columbus OH 43210.}
  ~\\
May, 1999
\end{flushright}

\vfill

{\bf ABSTRACT:} \break

The FKG theorem says that the POSITIVE LATTICE CONDITION,
an easily checkable hypothesis which holds for many natural families
of events, implies POSITIVE ASSOCIATION, a very useful property.  
Thus there is a natural and useful theory of positively dependent 
events.  There is, as yet, no corresponding theory of negatively
dependent events.  There is, however, a need for such a
theory.  This paper, unfortunately, contains no substantial theorems.  
Its purpose is to present examples that motivate a need for such a theory, 
give plausibility arguments for the existence of such a theory, outline a 
few possible directions such a theory might take, and state a number of 
specific conjectures which pertain to the examples and to a wish list of
theorems.  

\vfill



\end{titlepage}

\noindent{\bf Philosophy:} \\

The questions in this paper are motivated by several independent problems
in combinatorial probability, stochastic processes and statistical
mechanics.  For each of these problems, it seems that progress will require
(and engender) better understanding of what it means for a collection
of random variables to be ``repelling'' or mutually negatively dependent.
The temptation is to try to copy the theory of positively dependent 
random variables, since the FKG theorem and its offshoots give this 
theory a powerful footing from which to prove correlation inequalities, 
limit theorems and so on.  Perhaps it is folly: no definition of 
mutual negative dependence has proved one tenth as useful as the
lattice condition for positively dependent variables.  The purpose
of this paper is to lay the groundwork for whatever progress is possible
in this area.  The main goal is to state some conjectured implications 
which would bridge the gap between easily verifiable conditions and useful 
conclusions.  A second purpose is to collect together examples and
counterexamples that will be useful in forming hypotheses, and a third is
to update previous surveys by collecting the relevant known results and
adding a few more.  The scope of this paper is limited to binary-valued
random variables, in the hope that eliminating the metric and order 
properties of the real numbers in favor of the two point set $\{ 0 , 1 \}$ 
will better reveal what is essential to the questions at hand.

\section{Statement of the problem and some motivation}

\subsection{Definition of positive and negative association}

Let $\B_n$ be the Boolean lattice containing $2^n$ elements, each element
being thought of as a sequence of zeros and ones of length $n$, or as 
function from $\{ 1 , \ldots, n \}$ to $\{ 0 , 1 \}$, or as a subset
of $\{ 1 , \ldots , n \}$.  Let $\mu$ be a nonnegative function on the
lattice with $\sum_{x \in B_n} \mu (x) = 1$.  Then $\mu$ is a probability
measure on $\B_n$ and each coordinate function is a binary random variable,
denoted $X_j$, $j = 1 , \ldots , n$.  Sometimes we replace the base set
$\{ 1 , \ldots , n \}$ by a different index set arising naturally in an application,
such as the set of edges of a graph.

In order to make an analogy, we review the facts about positive
dependence.  The measure $\mu$ is said to be {\em positively associated} 
(c.f.\ Esary, Proschan and Walkup (1967)) if
\begin{equation} \label{eq PA}
\int fg \, d\mu \geq \int f \, d\mu \; \int g \, d\mu
\end{equation}
for every pair of increasing functions $f$ and $g$ on $\B_n$.  This is a 
strong correlation inequality from which many others may be
derived, and from which distributional limit theorems also follow;
see Newman (1980).  Positive association is implied by the following
local (and therefore often more checkable) {\em positive lattice condition} 
(Fortuin, Kastelyn and Ginibre (1971); see also Ahlswede and Daykin
(1979) for a more general proof):
\begin{th}[FKG] \label{th FKG}
If the following condition holds then $\mu$ is positively associated.
\begin{equation} \label{eq pos lattice}
\mu (x \vee y) \mu (x \wedge y) \geq \mu (x) \mu (y) .
\end{equation}
\end{th}
In fact, one only needs to check this in the case where $x$ and $y$
each cover $x \wedge y$ (an element $u$ covers an element $v$ if
$u > v$ and if $u \geq w \geq v$ implies $w \in \{ u , v \}$). 
This immediately allows verification of positive association for
basic examples such as the ferromagnetic Ising model, certain urn models, 
and, in the continuous case, multivariate normals, gammas, and many more
distributions.  Furthermore, the class of measures satisfying the
lattice condition~(\ref{eq pos lattice})
is easily seen to be closed under Cartesian products,
pointwise products, and, most importantly, under integrating out
any of the variables (i.e., any projection of $\mu$ onto the space
$\{ 0 , 1 \}^E$ for $E \subseteq \{ 1 , \ldots , n \}$ will also 
satisfy~(\ref{eq pos lattice})).

Negative dependence, by contrast, is not nearly as robust.  First,
since a random variable is always positively correlated with itself, one
cannot expect all monotone functions to be negatively correlated.  The
usual definition of {\em negative association} of a measure $\mu$
(c.f.\ Joag-Dev and Proschan (1983)) is that 
\begin{equation} \label{eq NA}
\int fg \, d\mu \leq \int f \, d\mu \int g \, d\mu
\end{equation}
for increasing functions $f$ and $g$,
provided that $f$ depends only on a subset $A$ of the $n$ variables and
$g$ depends only on a subset disjoint from $A$.  Secondly, whereas in
the positive case one may have $\E X_i X_j$ significantly greater than
$\E X_i \E X_j$ for many $i , j$, in the negative case the inequality
$\sum_{i,j} \Cov X_i X_j \geq 0$ prevents the typical term
$\Cov X_i X_j$ from having a significantly negative value.  Thirdly, 
the {\em negative lattice condition}, namely~(\ref{eq pos lattice}) 
with the inequality reversed, is not closed under projections.  
Thus one cannot expect it to imply negative association and indeed 
it does not.  

Contrasting the definitions of positive and negative association 
shows that the inequality~(\ref{eq PA}) comes from two sources.
The first is from autocorrelation when $f$ and $g$ depend on
the same variable in the same direction; thus for independent random
variables, strict inequality in~(\ref{eq PA}) occurs if $f$ and $g$
both depend on a common variable.  The second is from
positive interdependence of the variables which contributes even when 
$f$ and $g$ depend on disjoint subsets.  This leads immediately to
a question on positive association which, while not directly pertaining
to the subject of negative dependence, might shed light on how to 
disentangle inter- and auto-correlation.  
\begin{qq}
If one assumes~(\ref{eq PA}) only for $f$ and $g$ depending on 
disjoint subests of the variables, does the inequality follow 
for all increasing $f$ and $g$?  
\end{qq}
This elementary question has not, as far as I know, been posed or
answered in print.  

The reverse-inequality analogue of~(\ref{eq PA}) for product measures
is the van den Berg-Kesten-Reimer inequality:
\begin{equation} \label{eq BKR}
\mu (A \Box B) \leq \mu (A) \mu (B)
\end{equation}
Here $A \Box B$ is the event that $A$ and $B$ happen for ``disjoint
reasons'': $\omega \in A \Box B$ if there are disjoint subsets $S (\omega)$ 
and $T (\omega)$ of $\{ 1 , \ldots , n \}$ such that $A$ contains the 
set of all configurations agreeing with $\omega$ on $S$ and $B$ contains
the set of all configurations agreeing with $\omega$ on $T$.  This leads
to a different but also somewhat natural definition of negative association,
denoted here BKRNA (Berg-Kesten-Reimer negative association): a
measure $\mu$ has the BKRNA property if~(\ref{eq BKR}) holds for all holds
for all sets $A$ and $B$.

BKRNA has some claim to being ``the negative version'' of positive
association, since instead of reversing the inequality in~(\ref{eq PA})
and then restricting $f$ and $g$, we choose a different inequality to
reverse which holds in the independent case for all $f$ and $g$.  The
BKRNA property has been discussed in the literature, but has not
been fruitful.  This may be due to the fact that even in the independent
case, where the proof of~(\ref{eq PA}) has been known for 40 years
(see Harris 1960), the inequality~(\ref{eq BKR}) turned out to be 
quite hard to prove.  A proof when $A$ and $B$ are both up-sets 
(see definition next paragraph) was
given in van den Berg and Kesten (1985), generalized to the case where
$A$ and $B$ had the next level of complexity (up-set intersect down-set)
by van den Berg and Fiebig (1987), and then proved in complete generality
by Reimer in a manuscript yet to be published.  In view of this difficulty,
it seems unlikely that proving~(\ref{eq BKR}) for some interesting
non-product measure $\mu$ will be possible, let alone be the easiest 
way to establish a desired property of $\mu$.  Consequently, the remainder
of the paper deals with classical negative association, where we restrict
the test functions $f$ and $g$ instead of changing the binary set operation.

\subsection{Stochastic increase and decrease}

The notions of stochastic domination and stochastic increase and decrease
are useful when defining positive and negative dependence properties,
so we review them here.  Let $\mu$ and $\nu$ be measures on a 
partially ordered set, $S$.  An event $A \subseteq S$ is said to be
{\em upwardly closed} (or an {\em up-set}) 
if $x \in A$ and $y \geq x$ implies $y \in A$.
Often $S = \B_n$, the Boolean lattice of rank $n$, in which case this
is the same as $A$ being an increasing function of the coordinates.  
We say that $\mu$ {\em stochastically dominates} $\nu$ (written $\mu 
\stoch \nu$) if $\mu (A) \geq \nu (A)$ for every upwardly closed 
event $A$.  The condition $\mu_1 \stoch \mu_2 \stoch \cdots \mu_n$
is well known to be equivalent to the existence of a random sequence
$(X_1 , \ldots , X_n)$ such that $X_j \diseq \mu_j$ for each $j$ and 
$X_j \geq X_k$ for $1 \leq j \leq k \leq n$ (see e.g., Fill
and Michuda 1998).  We say that the random variable $X$ is 
{\em stochastically increasing} in the random variable $Y$ if the 
conditional distribution of $X$ given $Y = y_1$
stochastically dominates the conditional distribution of $X$ given
$Y = y_2$ whenever $y_1 \geq y_2$.  The notation $X \uparrow Y$ will
denote this relation, which is not in general symmetric.  Similarly,
$X$ is stochastically decreasing in $Y$ (denoted $X \downarrow Y$)
if one has $(X \| Y = y_1) \predeq (X \| Y = y_2)$ whenever $y_1 \geq y_2$.  
A convention in use throughout this paper is that terms involving inequalities 
are meant in the weak sense, so that for example ``decreasing'' means 
non-increasing and ``positively correlated'' means non-negatively correlated.

The relation $X \uparrow Y$ is not in general symmetric, but implies
$Y \uparrow X$ is a certain case, as given in the following proposition.
\begin{pr} \label{pr stoch sym}
Let $X$ be a $\{ 0 , 1 \}$-valued random variable and $Y$ take values
in any totally ordered set.  If $X \uparrow Y$ then $Y \uparrow X$.  
\end{pr}

\noindent{\sc Proof:}  Choose $t$ in the range of $Y$. Since 
$\P (X=1 \| Y)$ is inceasing in $Y$, it follows that 
$$\P (X=1 \| Y \leq t) \leq \sup_{s \leq t} \P (X=1 \| Y=s) \leq
   \inf_{s > t} \P (X = 1 \| Y=s) \leq \P (X=1 \| Y > t) .$$
Thus $X$ and $\one_{Y > t}$ are positively correlated and
$\P (Y > t \| X = 1) \geq \P (Y > t \| X = 0)$.  This holding
for all $t$ is equivalent to $Y \uparrow X$.   $\Cox$

A counterexample to the converse is given by the following 
probabilities, where the $(i,j)$-cell is the probability of
$(X,Y) = (i,j)$.  \\[2ex]
\hspace{1in}
\begin{tabular} {c|c|c|c|c|}
  & 1 & 2 & 3 & 4 \\ \hline
0 & 9/40 & 4/40 & 6/40 & 1/40 \\ \hline
1 & 1/40 & 6/40 & 4/40 & 9/40 \\ \hline
\end{tabular}

\subsection{Motivating examples}

The property of negative association is reasonably useful but hard
to verify.  The next subsection builds the case for ``reasonably useful''
by cataloging some consequences that would hold if negative dependence
could be established in some cases where it is conjectured.  In the present
subsection, we list some examples of systems which are known or
believed to have the negative association property.  The examples that
are conjectured motivate us to develop techniques for proving that measures
have negative dependence properties.  The point of including examples of
measures already known to be negatively associated is that we can use 
them to study properties of negative association, which will help us refine 
our conjectures about the consequences of negative association.  
As seen in Section~\ref{ss FM} below, knowledge of the characteristics 
of negatively associated variables will be helpful in proving criteria
for negative association.  

1.  The uniform random spanning tree.  Let $G$ be a finite connected graph,
and let $T$ be a random spanning tree (i.e. a maximal acyclic set of 
edges of $G$) chosen uniformly from among all spanning trees of $G$. 
It is easy to prove that the indicator functions $\{ X_e \}$ of the 
events that $e \in T$ have the following property: for any edges 
$e$ and $f$, $X_e$ and $X_f$ are negatively correlated.  Feder and
Mihail (1992) have shown that in fact this collection is negatively 
associated.  As we will see later, one concrete consequence of this is
that the conditional measures given $e \in T$ and $e \notin T$ may be
coupled to agree except that the latter has precisely one more edge 
elsewhere.  

A natural generalization is to consider weighted
spanning trees.  Let $W : E(G) \rightarrow \R^+$ be a function assigning
positive weights to the edges of $G$.  Define the weight $W(T)$ of a tree $T$ 
to be the product $\prod_{e \in T} W(e)$ of weights of edges in $T$.  The
probability measure $\mu$ on $\{ 0 , 1 \}^{E(G)}$ concentrated on 
spanning trees whose weights $\mu (T)$ are proportional to $W(T)$
is called the weighted spanning tree measure.  Everything known about
the uniform spanning tree also holds for the weighted spanning tree; in
fact a rational edge weight of $r/s$ may be simulated in the uniform
spanning tree setting by replacing the edge $e$ by $r$ parallel paths
of length $s$ each.  

2.  Simple exclusion.  Let $G$ be a finite graph, let $\eta_0$ be 
a function from $V(G)$ to $\{ 0 , 1 \}$, and let $\xi_t$ be the
trajectory of a simple exclusion process starting from $\xi_0 = \eta_0$.
The simple exclusion process is the Markov chain described as follows.
For each edge $e$ independently, at times of a rate 1 Poisson process,
the values of $\eta$ at the two endpoints of $e$ are switched.  This is
thought of as a particle moving across the edge but only if the opposite
site is vacant.  Fix $t$ and let $X_v = \xi_t (v)$ be the indicator function 
of the occupation of the vertex $v$ at time $t$.  It is known (Liggett
1977) that
\begin{equation} \label{eq cyl}
\E \left [ \prod_{v \in S} X_v \right ]  \leq \prod_{v \in S} \E X_v \, 
\end{equation}
for any subset $S$ of the vertices of $G$.  Are the variables $X_v$ 
negatively associated?  The most natural generalization of simple exclusion
is to allow the Poisson processes on the different edges to have different
rates; the inequality~(\ref{eq cyl}) is known in this generality.  

3.  Random cluster model with $q < 1$.  Let $G$ be a finite graph.  
For any subset $\eta$ of the edges, viewed as a map 
$\eta : E(G) \rightarrow \{ 0 , 1 \}$, let $N (\eta)$ denote the 
number of connected components of the graph represented by $\eta$.  
Given parameters $p \in (0,1)$ and $q > 0$, define a measure 
$\mu = \mu_{p,q}$ on $\{ 0 , 1 \}^E$ by letting 
\begin{equation} \label{eq RC}
\mu (\eta) = C p^{\sum_e \eta (e)} (1-p)^{\sum_e 1-\eta (e)} 
   q^{N(\eta)} \, .
\end{equation}
Here $C$ is the normalizing constant 
$$ C = \left [ \sum_{\eta : E(G) \rightarrow \{ 0 , 1 \}} p^{\sum_e \eta (e)} 
   (1-p)^{\sum_e 1-\eta (e)} q^{N(\eta)} \right ]^{-1}  .  $$
When $q > 1$, the variables $X_e := \eta (e)$ are easily seen to be
positively associated by checking the positive lattice condition and
applying the FKG Theorem.  When $q < 1$, the negative lattice condition
holds, but aside from this little is known about the extent of 
negative dependence.  Negative association and BKRNA are both conjectured
to hold, but it is not even known whether the variables $X_e := \eta (e)$
are pairwise negatively correlated under $\mu$.  The random 
cluster model has the uniform spanning tree model as a limit as $p,q$ 
and $p/q$ go to zero (see H\"aggstr\"om 1995); thus negative 
association in the RC model would in a way generalize what is known for 
spanning trees.  The RC model may be generalized by letting the factor $p$ 
vary from edge to edge.  Thus one has a function $p : E(G) \rightarrow (0,1)$
and the term $p^{\sum \eta (e)} (1-p)^{\sum 1 - \eta (e)}$ is replaced
by the more general $\prod_e p(e)^{\eta (e)} (1-p(e))^{1 - \eta(e)}$.

4.  Occupation of competing urns.  Let $n$ urns have $k$ balls dropped
in them, where the locations of the balls are IID chosen from some
distribution.  Let $X_i$ be the event that urn number $i$ is non-empty.
It is proved in Section~\ref{ss urns} that these events are negatively 
associated.  Dubhashi and Ranjan (1998) consider this example at length
and show negative association of the occupation numbers of the bins
(numbers of bals in each bin).  From this follows negative association
of the indicators of exceeding any prescribed threshholds $a_i$ in bin $i$.
Occupation numbers of urns under various probability schemes have appeared
many places.  Instead of multinomial probabilities, one can postulate
indistinguishability of urns or balls and arrive at Bose-Einstein or
other statistics.  Negative association seems only to arise in the
multinomial models, where Mallows (1968) was one of the first to
observe negative dependence.  

\subsection{Consequences of positive and negative association}

One use that is reasonably general is that of classifying
infinite volume limits of Gibbs measures.  The prototypical
example is the ferromagnetic Ising model.  The ferromagnetic 
Ising measure on a finite box $G$ with boundary $B$ and boundary
condition $\eta : B \to \{-1 , 1\}$ is a measure on spin configurations
$\xi : G \to \{ -1 , 1 \}$ proportional to
$$\exp \left ( \beta \left (\sum_{x,y \in G} \xi (x) \xi (y) + 
   \sum_{x \in G , y \in B} \xi (x) \eta (y) \right ) \right ) .$$
The spin variables $\{ \xi (x) : x \in G \}$ are positively associated and
stochastically increasing in $\{ \eta (y) : y \in B \}$,
from which it follows that there are a stochastically greatest and
least infinite volume limit, corresponding to plus and minus boundary
conditions respectively.  Thus there is non-uniqueness of the Gibbs
state if and only if the plus and minus states differ.  

Another example of this is the uniform spanning tree, which is almost
Gibbsian except that some configurations have infinite energy (are
forbidden).  Let $\mu_n^{(A)}$ be the uniform spanning tree measure
on the finite subcube of the $d$ dimensional integer lattice centered
at the origin with semi-diameter $n$.  The $A$ refers to a specification
of boundary conditions, i.e., of a partition of the vertices of the boundary
of the $n$-cube into components, so that the sample tree is uniform over
all spanning forests of the cube that become trees if each component
of $A$ is shrunk to a point.  Pemantle (1991) shows that the measures 
$\mu_n^{(A_n)}$ converge weakly to a measure $\mu$ in the case where $A_n$ 
is the discrete partition, and uses electrical network theory to show that 
this same limit holds for any $A_n$.  With the negative association
result of Feder and Mihail (1992) it is easy to see this directly as
follows.  Iterating the stochastic relation between the conditional
measures given $e \in T$ and given $e \notin T$ shows that
$\mu_n^{(A)} \predeq \mu_n^{(A')}$ whenever $A'$ refines $A$.  Thus
the measures $\mu_n^{(A)}$ are stochastically sandwiched between the
measures induced by ``free'' and ``wired'' boundary conditions (where
$A$ is repectively discrete or a single component); thus the set of
limits is sandwiched between a maximal and minimal limit measure;
both must have the same one-dimensional marginals (by stationarity)
and hence must coincide.  

Negative association has the further consequence that the 
uniform spanning tree measure is Very Weak Bernoulli.  
Briefly, this means that the conditional measures inside a large box 
given two independent realizations of the boundary can be coupled
so as to make the expected proportion of disagreements arbitrarily low.
To see that the Uniform Spanning Tree is VWB, note that the number of edges in 
a spanning tree is determined by the boundary conditions, so that free
boundary conditions will always yield precisely $|\partial B| - 1$ more
edges than wired boundary conditions, where $\partial B$ denotes the
set of vertices in the boundary of a set $B$.  Given two boundary
conditions $A_1$ and $A_2$, we can construct a triple $(T_1 , T_* , T_2)$
such that $T_1$ is chosen from the measure with boundary conditions 
$A_1$, $T_2$ from boundary conditions $A_2$, and $T_*$ from free
boundary conditions, and so that $T_*$ contains $T_1$ (construct
$(T_1, T_*)$ from the coupling witnessing $T_1 \predeq T_*$ and then
construct $T_2$ given $T_*$ from a coupling witnessing $T_2 \predeq T_*$).
Then $T_1$ and $T_2$ differ in fewer than $2 |\partial B|$ places.
Question: is there a simultaneous coupling of all boundary conditions
such that the configuration with boundary condition $A$ is a subset of
the configuration with boundary condition $A'$ whenever $A'$ refines $A$?
For the reason why this does not immediately follow from stochastic
monotonicity in the boundary conditions, see Fill and Machida (1998).  

Positive and negative association may be used to obtain information on
the distribution of functionals such as $\sum_e X(e)$.  Newman 
(1980, 1984) shows that under either a positive or negative dependence
assumption, of strength between cylinder dependence and full association,
the joint characteristic function of the variables $\{ X_e \}$ is well
approximated by the product of individual characteristic functions.
This allows him to obtain central limit theorems for stationary
sequences of associated variables.  In the positive association case 
one needs to assume summable covariances, whereas in the negative
case one gets this for free.  It is logical to ask what information may
be obtained from negative association without passing to a limit.  For
example, since one has a CLT or triangular array theorem in the
independent case, can one prove that negatively associated events
are at least as tightly clustered as independent events?  
Section~\ref{ss rank seq} discusses some conjectures along 
these lines.  Here is a specific application of these conjectures.

Consider simple exclusion on the one-dimensional integer lattice, with
initial configuration given by $X_v = 1$ for $v \leq 0$ and $X_v = 0$
for $v > 0$.  What can one say about the number $N_t := \sum_{v > 0}
\eta_t (v)$ of occupied sites to the right of the origin at time $t$? 
The mean $\E N_t$ is easy to compute, and an upper bound of
$O(t^{1/2})$ on the variance has been obtained by several people.
While this shows that $(N_t - \E N_t) / t^{1/4}$ is tight, it is 
a far cry from a limit theorem.  It would be nice to be able to obtain
a central limit theorem, or, in lieu of that, Gaussian bounds on the
tails of $N_t$.  The conjectured chain of implications is: 
first, the exclusion model is negatively associated; second, negatively 
associated measures have sub-Gaussian tails.  Negative assocation
is known [Dubhashi and Ranjan (1998), Proposition~7] to imply the 
Chernoff-Hoeffding tail bounds; see conjectures~(\ref{conj always ULC}) 
and~(\ref{conj ULC examples}) below for other possible consequences
of negative association.    

\subsection{Feder and Mihail's proof} \label{ss FM}

Feder and Mihail (1992) prove that a uniform random base for a 
{\em balanced matroid}, of which the uniform spanning tree measure
is a special case, has the negative association property\footnote{This is
false for general matroids; see Seymour and Welsh (1975).}.  They use 
induction on the size of the edge set $E$, with the specific nature of 
the measure entering through only two properties, $(i)$ and $(ii)$.  The
logical form of the proof is as follows.  Choose an edge $e$ appropriately
and show that property~$(ii)$ holds for $(\mu \| e)$.  This together with
property~$(i)$ for $\mu$ and the induction hypothesis then imply that 
$\mu$ is negatively associated.  

This argument provides further motivation for deriving consequences of
negative association.  If we can prove, for example, that negative
association implies property$~(ii)$, then the step where we verify 
property~$(ii)$ drops out (by induction!) and the entire argument may be
carried out using only property~$(i)$.  Proving something weaker 
than~$(ii)$ for negatively associated measures still reduces the
work to proving~$(ii)$ from this property.  We make this all concrete by
defining the properties and stating the above as a theorem.  

Let $S$ be a class of measures on Boolean algebras which is
closed under conditioning on some of the coordinate values.  An 
example of such a measure is the uniform or weighted spanning tree
measure or the random cluster measure.  
\begin{quote}
Property $(i)$ {\bf pairwise negative correlation}: each $\mu \in S$
makes each pair of distinct $X_e$ and $X_f$ negatively correlated.

Property $(ii)$ {\bf some edge correlates with each up-set}: for each
$\mu \in S$ and increasing event $A$ there is an edge $e$ with 
$\mu (X_e \one_A) \geq \mu (X_e) \mu(A)$.  
\end{quote}
\begin{th} \label{th FM argument}
Let $\S$ be a class of measures closed under conditioning and under
projection (i.e., forgetting some of the variables) and suppose 
all measures in this class have pairwise negative correlations.
Then property~$(ii)$ for $S$ (implied for example by 
Conjecture~\ref{conj rank cover} below) implies that every measure 
in $\S$ is negatively associated.  
\end{th}

\noindent{\sc Proof of theorem:}  Pick $\mu$ in $S$ and induct on
the rank $n$ of the lattice on which $\mu$ is a measure.  When $n=1$ 
the statement is trivial.  Now assume the conclusion for all measures
in $S$ on lattices of size less than $n$.  The remainder of the proof 
copies the Feder-Mihail argument.  For brevity, we show that $A$ and $B$ 
are negatively correlated when $B = X_e$ and $A$ is an arbitrary up-set 
not depending on the variable $X_e$.  

If $\P (X_e = X_f = 1) = 0$ for all $f \neq e$ the induction step is
trivial, so assume not.  By property~$(ii)$ for $(\mu \| e)$ there 
is some $f \neq e$ for which
\begin{equation} \label{eq find edge}
\mu (A \| X_e = X_f = 1) \geq \mu (A \| X_e = 1)\, .
\end{equation}
Now write
\begin{eqnarray*}
&& \mu (A \| X_e = 1) = \\
&&~~~\mu (X_f = 1 \| X_e = 1) \mu (A \| X_e = X_f = 1) + 
   \mu (X_f = 0 \| X_e = 1) \mu (A \| X_e = 1 , X_f = 0) \\[2ex]
&&\mu (A \| X_e = 0) = \\
&&~~~\mu (X_f = 1 \| X_e = 0) \mu (A \| X_e = 0 , 
   X_f = 1) + \mu (X_f = 0 \| X_e = 0) \mu (A \| X_e = 0 , X_f = 0) \, .
\end{eqnarray*}
Comparing terms on the right-hand sides, we see that 
\begin{quote}
$(i)$~$\mu (X_f = 1 \| X_e = 1) \leq \mu (X_f = 1 \| X_e = 0)$ by the 
assumption that measures in $\S$ have pairwise negative correlations; \\

$(ii)$~$\mu (A \| X_e = X_f = 1) \leq \mu (A \| X_e = 0 ,  X_f = 1)$
since the conditional law $(\mu \| X_f = 1)$ is assumed by induction to
be negatively associated and hence $A$ and $X_e$ are negatively correlated
given $X_f = 1$; \\

$(iii)$~$\mu (A \| X_e = 1 , X_f = 0) \leq 
\mu (A \| X_e = 0 , X_f = 0)$ by the induction hypothesis this time applied
to $(\mu \| X_f = 0)$; \\

$(iv)$~$\mu ( A \| X_e = X_f = 1) \geq 
\mu ( A \| X_e = 1 , X_f = 0)$ by the choice of $f$.  
\end{quote}
These four imply that the left-hand sides are comparable: 
$\mu (A \| X_e = 1) \leq \mu (A \| X_e = 0)$.  This completes the 
induction in the special case where one of the two upwardly closed events 
is a simple event, $\{ X_e = 1 \}$.  The case of a general upwardly closed 
event is similar (see the Exercise 6.10 in Lyons and Peres 1999).
$\Cox$

\section{Properties and implications}

\subsection{Obtaining measures from other measures} \label{ss closure}

Before discussing negative dependence properties of various
strengths, we consider ways of obtaining a measure $\mu'$
from a given measure $\mu$ in such a way as to preserve
any known or conjectured negative dependence properties.  
The reason for discussing these beforehand is to lend perspective
to some of the definitions: if the property is not closed under the
$\mu \mapsto \mu'$, either by definition or by some argument,
then perhaps it is not such a natural property.  In the foregoing,
we fix a finite set $E$ and a probability measure $\mu$ on 
the space $\{ 0 , 1 \}^E$.

1.  Projection.  Given $E' \subseteq E$, let $\mu'$ be the projection
of $\mu$ onto $\{ 0 , 1 \}^{E'}$.  This corresponds to integrating 
out (i.e., forgetting) the variables in $E \setminus E'$.  
Clearly any natural negative dependence property is closed under projection.

2.  Conditioning.  Given $A \subseteq E$ and $\eta \in \{ 0 , 1 \}^A$,
consider the conditional distribution $(\mu | X_e = \eta (e) \mbox{ for }
e \in A)$.  It is reasonable to expect these sections of the measure $\mu$
to be negatively dependent if $\mu$ is.  Several of the motivating examples,
namely spanning trees, RC model and the Ising model, are classes of
measures closed under conditioning.  Note that we are not allowing
conditioning on a set larger than a single atom.  To ask that the projection
of $\mu$ onto $\{ 0 , 1 \}^{E \setminus A}$ be negatively dependent, 
conditioned on the event $< X_e : e \in A > \in S$ for arbitrary $S$ is
significantly stronger.  

3.  Products.  If $\mu_1$ and $\mu_2$ are negatively dependent, then
clearly $\mu_1 \times \mu_2$ should be.  

4.  Relabeling.  The measure $\mu'$ defined by $\mu' \{ X_e = \eta (e) :
e \in E \} = \mu \{ X_e = \eta (\pi (e)) : e \in E \}$, where $\pi$
is some permutation of $E$, is of course just a relabeling of $\mu$.

5.  Extends the concept of negative correlation.  When $|E| = 2$, 
any reasonable definition reduces to negative correlation.  

6.  External field.  The name for this property is borrowed from
the Ising model.  Let $W : E \rightarrow \R^+$ be a non-negative
weighting function and let $\mu'$ be the reweighting of $\mu$ by $W$.
Specifically, let 
$$\mu' \{X_e = \eta (e): e \in E \} = C \prod_{e \in E} W(e)^{\eta (e)}
   \mu \{X_e = \eta (e): e \in E \} \, ,$$
where $C$ is a normalizing constant.  This corresponds to making
a particular value for each edge more or less likely, without introducing
any further interaction between the edges.  For example if $W(e) \neq 1$
for a unique $e$, then the probability of $\{ X_e = 1 \}$ is altered, but
the conditional distributions of $(\mu | X_e)$ are unaltered.  Many of the
classes of measures which motivate our study are closed under imposition
of an external field.  For spanning trees or for the RC model, this 
corresponds to the weighted case; for the Ising model it corresponds 
to an external field.  Closure under external fields may seem far from a
natural condition for models that are not thermodynamic ensembles, but
this may be more natural than it seems.  First, if one believes in
closure under conditioning, then this is the canonical interpolation between
conditioning on $X_e = 1$ and conditioning on $X_e = 0$.  Secondly, 
Karlin and Rinott in 1980 had already proposed a property they call
$\neg$ which is essentially the negative lattice condition
plus closure under projection and external fields (see the discussion
preceding Conjecture~\ref{conj horizontal}).  

\subsection{Negative dependence properties and their relations}
\label{ss NDP}

We recall the definition of negative association:
\begin{defn} 
$\{ X_e : e \in E \}$ are negatively associated (NA) if for every $A 
\subseteq E$ and every pair of bounded increasing functions
$f: \{ 0 , 1 \}^A \rightarrow \R$ and $g : \{ 0 , 1 \}^{E \setminus A} 
\rightarrow \R$, $\E fg \leq \E f \E g$.  
\end{defn}
Unfortunately, this property is not closed under conditioning or external 
fields (see Example 2 below).
This may be an indication that these two closures are not so natural
after all, but on the other hand it makes sense, at least for 
closure under conditioning, to make a new definition:
\begin{defn}
The measure $\mu$ is conditionally negatively associated (CNA) if 
each measure $\mu'$ gotten from $\mu$ by conditioning on some (or none)
of the values of the variables is negatively associated. 
\end{defn}
Since the operation of conditioning is easy to understand in many of
our motivating examples, this extension should not prove to unwieldy.

The weakest possible negative dependence property is pairwise negative
correlation: \\ 
$\mu (X_e X_f) \leq \mu (X_e) \mu (X_f)$.  For real-valued
random variables, there is a stronger pairwise property, called
{\em negative quadrant dependence} (NQD) in Newman (1984),  
after Lehman (1966).  Say that $X$ and $Y$ are NQD if 
$$\P (X \geq a , Y \geq b) \leq \P (X \geq a) \P (Y \geq b)$$
for all $a$ and $b$.  For binary-valued random variables, this
reduces to simple correlation.  A stronger property, called
{\em negative regression dependence} (in analogy with positive 
regression dependence c.f.\ Esary, Proschan and Walkup 1967),
is defined by requiring the conditional distribution of $X$ given
$Y$ to be stochastically decreasing in $Y$: $\P(X \geq t | Y = s)$
is decreasing in $s$ for each $t$.  For binary-valued variables this
again reduces to negative correlation.  When $X$ and $Y$ are vectors,
$X := < X_e : e \in A > , Y := <x_e : e \notin A>$, this would say
that the conditional joint distribution of $\{ X_e : e \in A \}$ 
given $\{ X_e : e \notin A \}$ should be stochastically decreasing
in the values conditioned on.  Thus we have a definition:
\begin{defn}
Say that the variables $\{ X_e : e \in E \}$ are jointly negative
regression dependent (JNRD) if the vectors $<X_e : e \in A >$ and
$<X_e : e \notin A >$ are always negative regression dependent.  
Equivalently, require that for any increasing event $H$ measurable 
with respect to $\{ X_e : e \in A \}$, 
$\mu (H | x_e : e \notin A)$ is decreasing with respect to the
partial order on $\{ 0 , 1 \}^{A^c}$.  
\end{defn}

Unraveling the definitions, one sees that conditional negative 
association implies JNRD, since JNRD is simply CNA in the
special case where one has conditioned on $\{ X_e : e \in A^c \}
\setminus \{ f \}$ and then asks for $X_f$ to be negatively
correlated with $\one_H$ for any increasing event $H$ measurable
with respect to $\{ X_e : e \in A \}$.  

The negative lattice condition 
\begin{equation} \label{eq neg lattice}
\mu (x \vee y) \mu (x \wedge y) \leq \mu (x) \mu (y) .
\end{equation}
is closed under five of the six closure operations, but the missing
one, projection, is crucial.  This is what makes the negative version
of the FKG theorem fail.  Accordingly, 
\begin{defn}
Say that $\{ X_e : e \in E \}$ satisfy the hereditary negative lattice
condition (h-NLC) if every projection satisfies the negative lattice
condition.  
\end{defn}
It is easy to see that JNRD implies h-NLC, since h-NLC is the special
case where $A$ is a singleton.  

None of the three properties CNA, JNRD or the hereditary NLC are 
closed under imposition of an external field (see Example 1 below).
Projecting from index set $S$ to $S'$ and then imposing an external field 
(on $S'$) is the same as imposing an external field which is trivial on 
$S \setminus S'$ and then projecting to $S'$.  Thus any sequence of
projections and external fields may be written as one external field
followed by one projection.
One may define three stronger properties, CNA+, JNRD+ and h-NLC+, 
which are that the corresponding properties hold for
the given measure and for all measures obtained from the given measure
by imposition of an external field and a projection; these properties
are then by definition closed under external fields and projections.  
While these stronger properties are difficult to check directly, they
appear to hold for the motivating examples and are introduced
in the hope that they do in fact hold there and are strong enough
to be useful in inductive arguments such as the proof of 
Theorem~\ref{th FM argument}.  The property h-NLC+
is called $\neg$ by Karlin and Rinott (1980), according
to terminology they develop mainly for continuous random variables. 

The terminology introduced thus far can be summarized with a diagram 
of implications. \\[3ex]

\setlength{\unitlength}{2pt}
\begin{picture}(170,80)(20,-20)
\put(30,60){CNA+}
\put(80,60){JNRD+}
\put(130,60){h-NLC+}
\put(130,67){($\neg$)}
\put(30,30){CNA}
\put(80,30){JNRD}
\put(130,30){h-NLC}
\put(30,0){NA}
\put(49,62){\vector(1,0){27}}
\put(99,62){\vector(1,0){27}}
\put(49,32){\vector(1,0){27}}
\put(99,32){\vector(1,0){27}}
\put(34,56){\vector(0,-1){20}}
\put(84,56){\vector(0,-1){20}}
\put(134,56){\vector(0,-1){20}}
\put(34,26){\vector(0,-1){20}}
\put(90,-10){Figure 1}
\end{picture}

\subsection{Conjectures, examples and counterexamples} \label{ss urns}

The vertical implications in Figure~1 are strict, as shown by the examples 
which follow in this section.  Whether the horizontal implications are
strict is an open question:
\begin{conj} \label{conj horizontal}
All three properties CNA+, JNRD+ and h-NLC+ are equivalent.
\end{conj}
Another immediate question is whether anything other than CNA is
strong enough to imply negative association.  
\begin{conj} \label{conj 2}
Strong version: h-NLC implies NA. Weak version: h-NLC+ implies NA.
\end{conj}

Examples showing the vertical implications are not equivalences are as
follows (verified by brute force). 

\noindent{\em Example 1:} Suppose $n = 3$, and the probabilities for
the various possible atoms are proportional to the following: 
\begin{eqnarray*}
\P(X_1 = 0 , X_2 = 0 , X_3 = 0) & = & 16 \\ 
\P(X_1 = 0 , X_2 = 0 , X_3 = 1) & = & 8 \\ 
\P(X_1 = 0 , X_2 = 1 , X_3 = 0) & = & 8 \\ 
\P(X_1 = 0 , X_2 = 1 , X_3 = 1) & = & 8 \\ 
\P(X_1 = 1 , X_2 = 0 , X_3 = 0) & = & 12 + \ee  \\ 
\P(X_1 = 1 , X_2 = 0 , X_3 = 1) & = &4 \\ 
\P(X_1 = 1 , X_2 = 1 , X_3 = 0) & = & 4 \\ 
\P(X_1 = 1 , X_2 = 1 , X_3 = 1) & = & 1  \, .
\end{eqnarray*}
When $0 \leq \ee \leq .8$ then this measure satisfies CNA and
hence JNRD and h-NLC.  However, when $\ee > 0$, then applying 
the external field $(\lambda , 1 , 1)$ for any positive   
$\lambda < \ee / (1 - \ee)$ yields a measure in which $X_2$
and $X_3$ are positively correlated, thus violating h-NLC
and hence JNRD and CNA.  This shows the first three vertical implications 
in Figure~1 are strict.

\noindent{\em Example 2:} Suppose $n = 3$, and the probabilities for
the various possible atoms are in the proportions:
\begin{eqnarray*}
\P(X_1 = 0 , X_2 = 0 , X_3 = 0) & = & 0 \\ 
\P(X_1 = 0 , X_2 = 0 , X_3 = 1) & = & 1 \\ 
\P(X_1 = 0 , X_2 = 1 , X_3 = 0) & = & 1 \\ 
\P(X_1 = 0 , X_2 = 1 , X_3 = 1) & = & 10 \ee  \\ 
\P(X_1 = 1 , X_2 = 0 , X_3 = 0) & = & 1  \\ 
\P(X_1 = 1 , X_2 = 0 , X_3 = 1) & = & 1  \\ 
\P(X_1 = 1 , X_2 = 1 , X_3 = 0) & = & 10 \ee  \\ 
\P(X_1 = 1 , X_2 = 1 , X_3 = 1) & = & \ee  \, .
\end{eqnarray*}
Here the negative lattice condition fails on the four atoms 
having $X_2 = 1$; thus CNA, JNRD and h-NLC (in fact NLC) all fail,
whereas the variables are in fact negatively associated.  Thus the
lowest vertical implication in Figure~1 is strict as well.

The following lemma will be useful on a number of occasions.  The easy
inductive proof is omitted.

\begin{lem} \label{lem Markov}
Let $Y_1 , \ldots , Y_n$ be random variables taking values in a partially
ordered set and suppose they have the Markov property, namely that $Y_1 ,
\ldots , Y_{k-1}$ are independent from $Y_{k+1} , \ldots , Y_n$
given $Y_k$.  Suppose also that each $Y_{k+1}$ is either 
stochastically increasing or decreasing in $Y_k$.  Then
$Y_n$ is either stochastically increasing in $Y_1$ or
stochastically decreasing in $Y_1$, according to whether
the number of indices $k$ for which $Y_{k+1}$ is decreasing in $Y_k$
is even or odd.   $\Cox$
\end{lem}

We conclude this subsection with a proof that the competing urn model 
of Example~4 is negatively associated.  The result with general
threshholds is proved in Dubhashi and Ranjan (1998), but the
proof given here is independent of that.  

\noindent{\sc Proof that the urn model is negatively associated:}
Fix $1 < r < n$ and let $A$ and $A'$ 
be up-events measurable with respect to $\{ X_i : i \leq r \}$
and $\{ X_i : i > r \}$ respectively.  Let $V$ and $V'$ be the 
total number of balls dropped into urns $i$ with $i \leq r$ and
$i > r$ respectively.
Letting $Y_1$ be the indicator function of $A$, $Y_4$ be
the indicator function of $A'$, $Y_2 = V$ and $Y_3 = V'$,
it is clear that $Y_1 , Y_2 , Y_3 , Y_4$ has the Markov
property.  I claim also that $A$ is stochastically increasing
in $V$ and $A'$ is stochastically increasing in $V'$.  By symmetry,
consider only $A$ and $V$.  Observe that conditional
on $V = m$, the draws are exchangeable in the usual sense
(definition below), so we may condition
on the first $m$ draws being those that went in urns $i \leq r$.
Then the distribution of balls given $V = m$ and the distribution
of balls given $V = m+1$ may be coupled so that the latter
is always the former plus an extra ball somewhere.  This
establishes the claim.  It is similarly easy to show that 
$V'$ is stochastically decreasing in $V$.  By Proposition~\ref{pr stoch sym},
$V$ is stochastically increasing in $A$.  Then the hypothesis of the 
above lemma is satisfied with stochastic increase for $k=1$ and $k=3$
and stochastic decrease when $k=2$; it follows that
$A'$ is stochastically decreasing in $A$ which
proves negative association.  
$\Cox$

\subsection{The exchangeable case and the rank sequence}
\label{ss rank seq}

The variables $\{ X_1 , \ldots , X_n \}$ are said to be {\em exchangeable}
if their joint distribution is invariant under permutation.  In the
case of binary-values random variables, this is the same as saying 
that $\mu \{ X_k = \eta (k) : 1 \leq k \leq n \}$ depends only on
$\sum_k \eta (k)$.  A fair amount of intuition may be gained from 
this special case.  The conjectured equivalences in Figure~1 are
proved in this case, but more importantly, new conjectures come to 
light that ought to hold in the general case as well.  

For a measure $\mu$ on $\B_n$, define the {\em rank sequence}
$\{ a_k : 0 \leq k \leq n \}$ by $a_k := \mu \{ \sum_{j=1}^n X_j = k \}$.
Thus $\{ a_k : 0 \leq k \leq n\}$ gives the total probabilities for the $n+1$ 
ranks of the Boolean lattice $\B_n$.  If the random variables $\{ X_j \}$
are exchangeable, then $\mu$ is completely characterized by its rank 
sequence, with $\mu \{ X_j = \eta (j)  : 1 \leq j \leq n \} = a_k / 
{n \choose k}$ for $k = \sum_j \eta (j)$.  In this case, the 
negative lattice condition~(\ref{eq neg lattice}) boils down to 
log-concavity of the sequence $\{ a_k / {n \choose k} \}$ 
(a positive sequence is said to be {\em log-concave} if $a_k^2 \geq 
a_{k-1} a_{k+1}$).  This motivates the following definition.
\begin{defn}
A finite sequence $\{ a_k : 0 \leq k \leq n \}$ is said to be 
Ultra-Log-Concave (ULC) if the nonzero terms of the sequence 
$\{ a_k / {n \choose k} \}$ form a log-concave sequence
and the indices of the nonzero terms form an interval.
\end{defn}

\noindent{\bf Convention:} From now on, to avoid trivialities,
we have included in the definition of log-concavity that the indices of 
the nonzero terms form an interval.  It will be useful later to note
that log-concavity is conserved by convolutions and pointwise
products.  

The significance of Ultra-Log-Concavity in the general case is still
conjectural, but in the exchangeable case it is given by the following
theorem whose proof appears at the end of the section.  
\begin{th} \label{th exchangeable}
Suppose that $\{X_j \}$ are exchangeable.  Then the six conditions 
CNA+, JNRD+, h-NLC+, CNA, JNRD and h-NLC (see Figure~1) 
are equivalent to Ultra-Log-Concavity of the rank sequence $\{ a_k \}$. 
This is trivially equivalent to the negative lattice 
condition,~(\ref{eq neg lattice}).  
\end{th}

Call the measure $\mu$ (not necessarily exchangeable) a ULC measure if 
its rank sequence is
ULC, and use the term ULC+ to denote a measure such that any measure
obtained from it by external fields and projections is ULC.  The
following conjectures, if true, imply a large role for the ULC 
property in the study of negative dependence.  They have been checked
only for lattices of rank up to~4.  

\begin{conj} \label{conj always ULC}
The strongest version of this conjecture is that any negatively associated
measure is ULC.  For a weaker version, replace the hypothesis of NA by
any of the other six stronger conditions in Figure~1.
\end{conj}
\begin{conj} \label{conj ULC examples}
In the RC model, the sum $\sum_{e \in S} X_e$ over any subset $S$ has 
a ULC rank sequence.  The same holds for the competing urns model.  In the
exclusion model, the total number of occupied sites in any set $S$ at any 
time $t$ has ULC rank sequence.
\end{conj} 

\noindent{\em Remark:} The ULC property for number of edges present
 from a given subset in a uniform (or weighted) random spanning tree 
is a subcase of the conjecture for the RC model.  For spanning trees,
this would sharpen a result of Stanley (1981) showing that the rank
sequence for a uniform random base of a unimodular matroid (of
which the uniform spanning tree is a special case) is log-concave.

Conjecture~\ref{conj always ULC} or the weaker~\ref{conj ULC examples}
would serve two purposes.  Firstly, the ULC property implies tail 
estimates on a distribution.  Secondly, Conjecture~\ref{conj always ULC} would
imply that that the ULC property is a necessary condition for
negative association, which helps to narrow and define our search
for the ``right'' negative dependence property.

The fact that ULC implies CNA+ {\em et al} in the exchangeable case
leads one to believe that ULC+ might be enough to imply negative
dependence in general:
\begin{conj} \label{conj ULC implies NA}
If $\mu$ is ULC+ then $\mu$ is CNA (hence CNA+) and in particular
$\mu$ is negatively associated. 
\end{conj}
Unlike the previous two, this conjecture is not particularly useful, 
since the hypothesis of ULC+ is hard to check.  It would, however, 
have philosophical value: supposing there to be a useful definition
of negative dependence still lurking out there, we have been approximating
it from the weak side, finding criteria that certainly hold for any 
such definition; the foregoing conjecture strengthens our previous
approximation by adding the property ULC+.  

A final philosophical observation belongs in this section.  
If Ultra-Log-Concavity is, as conjectured, a property of all 
negatively dependent measures, then the class of ULC sequences must 
be closed under convolution.  Indeed, if $\mu_1$ and $\mu_2$ 
are two exchangeable measures with ULC rank sequences,
then by Theorem~\ref{th exchangeable} they are negatively dependent 
in all senses we can imagine, so their product must be as well.  The
rank sequence for the product is the convolution of the rank sequences,
so unless even our understanding of the exchangeable case is nil,
the following conjecture must be true.  Embarrassingly, in the previously
circulated draft of this paper, there was no proof of the following
conjecture.  It has recently been proved by Liggett (1997).
\begin{conj}[Now proved by Liggett] \label{conj ULC convolve}
The convolution of two ULC sequences is ULC.  
\end{conj}

This section concludes with a proof of Theorem~\ref{th exchangeable}.
Begin with the following two lemmas.
\begin{lem} \label{lem stays ULC}
Let $\mu$ be an exchangeable measure with ULC rank sequence.  Suppose the 
measure $\mu'$ is obtained from $\mu$ by imposing an external field 
at coordinates $1 , \ldots k$ (i.e., $W(j) = 1$ for $j > k$) and
then projecting onto coordinates $r+1 , \ldots , n$ for some $r \geq k$.
Then $\mu'$ is exchangeable with ULC rank sequence.
\end{lem}

\noindent{\sc Proof:}  The exchangeability of $\mu'$ is clear.  To see
that $\mu'$ has ULC rank sequence, it suffices to consider the case 
$r = 1$.  [Reason: defining $\mu_j$ to be the measure gotten by imposing the
external field on the first $j$ coordinates and projecting onto the
last $n-j$ coordinates, one sees by induction on $j$ that $\mu_r
= \mu'$ will have the desired property].  So we assume without loss
of generality that $k = r = 1$.  

Let $\lambda$ denote $W(1)$.  Let $a_j$ (respectively $a_j'$) denote
the rank sequence for $\mu$ (respectively $\mu'$) and let $q_j$
(respectively $q_j'$) denote $a_j / {n \choose j}$ (respectively 
$a_j' / {n-1 \choose j}$).  Then 
$$q_j' = C (q_j + \lambda q_{j+1})\, , $$
where $C$ is the normalizing constant for the external field.  By 
assumption, $\{ q_j \}$ is log-concave, and hence for any $i < j$,
$q_i q_j \leq q_{i+1} q_{j-1}$.  The proof is now a simple calculation.
\begin{eqnarray*}
&&C^{-2} \left [ ( q_j')^2 - q_{j-1}' q_{j+1}' \right ] \\[2ex]
& = & q_j^2 + 2 \lambda q_j q_{j-1} + \lambda^2 q_{j-1}^2 
   - q_{j-1} q_{j+1} - \lambda q_{j-2} q_{j+1} - \lambda 
   q_{j-1} q_{j} - \lambda^2 q_{j-2} q_j \\[2ex]
& = & [ q_j^2 - q_{j+1} q_{j-1} ] + \lambda [ q_j q_{j-1} - 
   q_{j+1} q_{j-2} ] + \lambda^2 [ q_{j-1}^2 - q_j q_{j-2}] \, .
\end{eqnarray*}
This is the sum of three positive quantities, so it is positive,
proving log-concavity of $\{ q_j' \}$ which is equivalent to $\{ a_j' \}$
being ULC.   $\Cox$

\begin{lem} \label{lem SD for ex}
Let $\mu^*$ be a measure obtained from an exchangeable measure 
$\mu'$ with rank sequence $\{ a_k' \}$ by imposing an external field
$W$.  Let $Y_1$ and $Y_4$ be the respective indicator functions of 
$A$ and $A'$, events measurable with respect to disjoint sets $S$ and $S'$.  
Let $Y_2 = \sum_{e \in S} X_e$ and $Y_3 = \sum_{e \in S'} X_e$.  Then
the sequence $\{ Y_i \}$ is Markov.  Furthermore, the conditional laws 
$(\mu^* \| \sum X(e) = k)$ are stochastically increasing in $k$ and
the same holds for any projection of $\mu^*$ in place of $\mu^*$.
\end{lem}

\noindent{\sc Proof:}  Let $\mu', \mu^*, A, A' , S, S'$ and
$\{ Y_i \}$ be as in the hypotheses.  The probabilities
for $\mu^*$ are given as follows, with $C$ being a normalizing
constant as usual:
$$\mu^* \{ X_e = \eta (e), \mbox{ all } e \in S \} = C \prod_e W(e)^{\eta (e)} 
   {a_k' \over {n \choose k}} ,$$
where $k = \sum_e \eta (e)$.  From this, one gets the conditional 
probability
$$\mu^* ( X_e = \eta (e) : e \notin S' \| X_e = \eta (e) : e \in S') 
   = C' \prod_{e \notin S'} W(e)^{\eta (e)}  {a_k' \over {n \choose k}} .$$
This does not depend on the values of $\eta$ on $S'$ except through
$\sum_{e \in S'} \eta (e)$, which proves the Markov property.
For the stochastic increase, note that the conditional distribution
of $\mu^*$ given $\sum_e X(e)$ are the same as the law of independent
Bernoulli random variables with $\P (X(e) = 1) = W(e) / (1 + W(e))$,
conditioned on $\{ \sum_e X(e) = k \}$.  The same holds for any projection
of $\mu^*$.  There are elementary proofs that these laws increase 
stochastically in $k$, but in the context of this paper, the easiest
argument is to add an extra variable $X(e^*)$ and apply the 
Feder-Mihail result to the balanced matroid gotten by conditioning
on $\sum^*_e X(e) = k+1$ and to the conditional measures given
$X(e^*) = 0$ and $X(e^*) = 1$.    $\Cox$

\noindent{\sc Proof of Theorem}~\ref{th exchangeable}:  It is clear
that ULC is equivalent to the negative lattice condition and hence
is implied by h-NLC.  To show that ULC implies the other six conditions
we work up the ladder.  First, if $\mu$ is exchangeable and ULC, then
Lemma~\ref{lem stays ULC} shows that all projections of $\mu$ are as well, 
which means that the NLC holds hereditarily, giving h-NLC.  In fact,
the lemma is enough to give h-NLC+, since any $\mu^*$ obtained from
$\mu$ may be described (after re-ordering of coordinates) as some 
measure $\mu'$ as in the lemma, on which has been imposed an external field
(that is, any sequence of external fields and projections may be
written as an external field that affects only those indices not appearing
in the final measure, followed by a single projection, followed by an
external field); Lemma~\ref{lem stays ULC} implies $\mu'$ satisfies the 
negative lattice condition~(\ref{eq neg lattice}); this is invariant 
under external fields, so $\mu^*$ satisfies~(\ref{eq neg lattice}) as well.  

Next, we show that for any measure $\mu^*$ obtained from an exchangeable
measure $\mu$ by external fields and projections, JNRD implies CNA.
This will show that JNRD+ implies CNA+ as well as showing JNRD implies CNA.
To show this, let $\mu^*$ be such a measure.  Let $A$ and $A'$
be any up-events measurable with respect to disjoint sets of coordinates 
$S$ and $S'$.  Define a sequence of random variables $Y_1 , Y_2,
Y_3, Y_4$ by letting $Y_1$ be the indicator of $A$, letting $Y_4$ be the
indicator of $A'$, letting $Y_2 = \sum_{e \in S} X_e$, and letting
$Y_3 = \sum_{e \in S'} X_e$.  Apply Lemma~\ref{lem SD for ex} to
see that $\{ Y_i \}$ is Markov.  Lemma~\ref{lem Markov} finishes
the argument once we know that $Y_2$ is stochastically increasing
in $Y_1$, $Y_3$ is stochastically decreasing in $Y_2$, and $Y_4$
is stochastically increasing in $Y_3$.  Applying the last statement
of Lemma~\ref{lem SD for ex} to the projection of $\mu^*$ onto 
$\{ 0 , 1 \}^{S'}$, we see that the conditional joint law of 
$\{ X(e) : e \in S' \}$ given $\sum_{e \in S'} X(e) = k$ 
increases stochastically in $k$, which says precisely that
$Y_4$ is stochastically increasing in $Y_3$.  The same argument
with $S$ in place of $S'$ shows that $Y_1$ is stochastically
increasing in $Y_2$.  By Proposition~\ref{pr stoch sym},
$Y_2$ is stochastically increasing in $Y_1$.  Finally, to see that
$Y_3$ is stochastically decreasing in $Y_2$, write the conditional
distribution of $Y_3$ given $\{ Y_2 = k \}$ as an integral
$$\int {\rm Law} (Y_3 \| X(e) = \eta (e) : e \in S ) \, d \nu (\eta) ,$$
where $\nu$ is the mixing measure 
$$\nu \{ \eta \} = \mu^* (X(e) = \eta (e) : e \in S \| \sum_{e \in S} 
   X(e) = k) . $$
We have seen that $\nu$ is stochastically increasing in $k$.  By the
hypothesis that $\mu^*$ is JNRD, the integrand decreases stochastically
when $\eta$ increases in the natural partial order, and hence the
integral stochastically decreases in $k$.  This finishes the proof 
that JNRD implies CNA.  

It remains to show that h-NLC (respectively h-NLC+) implies
JNRD (respectively JNRD+).  The + case will be shown in 
Section~\ref{ss closed class} below, in the proof of 
Theorem~\ref{th closed class}, so we prove here only that 
ULC implies JNRD for exchangeable measures.  It suffices to 
show that the conditional distribution of $\sum_{e \neq f} X(e)$
given $X(f) = 0$ stochastically dominates the distribution of 
$\sum_{e \neq f} X(e)$ given $X(f) = 1$, since in the definition
of JNRD, comparing the conditional probabilities of any two
neighbors in the Boolean lattice $\{ 0 , 1 \}^{A^c}$ reduces
to comparing conditional probabilities given one value $X(f)$ 
after conditioning on all other values of $X(g) , g \in A^c$, and
such conditioning produces another exchangeable ULC measure.  
It further suffices to show that $\sum_{e \neq f} X_e$ is 
stochastically decreasing in $X_f$, since this is sufficient for
the distribution of $\{ X(e) : e \neq f \}$ given $X(f)$.  

Let $\{ a_j \}$ be the rank sequence for a ULC exchangeable measure
$\mu$, and let $\{ q_j \}$ be the sequence $\{ a_j / {n \choose j} \}$
as before.  Then 
$$\mu (\sum_{e \neq f} X_e = r \| X_f = 0) = {{n - 1 \choose r} q_r
   \over \mu (X_f = 0)} $$
and
$$\mu (\sum_{e \neq f} X_e = r \| X_f = 1) = {{n - 1 \choose r} q_{r+1}
   \over \mu (X_f = 1)}  \, .$$
Thus we need to show that for all $k < n$, 
$$\sum_{r=0}^k {{n-1 \choose r} q_{r+1} \over \mu (X(f) = 1)} \geq
  \sum_{r=0}^k {{n-1 \choose r} q_r \over \mu (X(f) = 0)} \, .$$
Cross-multiply and replace the quantities $\mu (X (f) = x)$ with
the sum over $s$ of $\mu (X(f) = x , \sum_{e \neq f} X(e) = s)$ to
transform this into
$$\sum_{r \leq k ; s \leq n-1} {n-1 \choose r} {n-1 \choose s} q_{r+1} q_s 
   \geq \sum_{r \leq k ; s \leq n-1} {n-1 \choose r} {n-1 \choose s} 
   q_r q_{s+1} .$$
Canceling terms appearing on both sides reduces the range of the
sum to $r \leq k < s$.  But for $r < s$, log-concavity of $\{ q_j \}$
implies that $q_{r+1} q_s \geq q_r q_{s+1}$, which establishes the
last inequality via term-by-term comparison and finishes the
proof that ULC implies JNRD.   $\Cox$

\section{Inductively defined classes of negatively dependent measures}

At this point it is worth examining the possibility that the many negative
dependence properties in our desiderata are not mutually satisfiable.
It is easy to see from the definition that the class of CNA+ measures 
is closed under products, projections and external fields, so we have 
at least one existence result:
\begin{quote}
Let $\S_0$ be the smallest class of measures containing all exchangeable 
ULC measures and which is closed under products, projections and external 
fields.  Then $\S_0$ is contained in the class of CNA+ measures.    $\Cox$
\end{quote}
Supposing there to exist a natural and useful class of ``negatively 
dependent measures'', it is contained in the class of CNA+ measures, 
and certainly contains the class $\S_0$.  This section aims 
to improve the latter bound which seems, intuitively to be further
 from the mark.  

\subsection{Further closure properties}

The class $\S_0$ is trivial, since products commute with external 
fields, and therefore $\S_0$ may be seen to contain only products 
of exchangeable ULC measures, on which have been imposed external fields.
We may enlarge the class $\S_0$ either by including more measures in the
base set or by increasing the number of closure operations in the inductive
step.  I will begin the discussion with a list of additional candidates for 
closure properties to those already listed in Section~\ref{ss closure}.  

7.  Symmetrization.  Given a measure $\mu$ on $\B_n$, let $\mu'$ 
be the exchangeable measure with $\mu' (\sum_j X_j = k) = 
\mu (\sum_j X_j = k)$.  In other words, $\mu' = (1 / n!) 
\sum_{\pi \in S_n} \mu \circ \pi$.  Since the measure $\mu'$
is exchangeable, we know criteria for $\mu'$ to be negatively
associated, and therefore closure under symmetrization boils
down to the Conjecture~\ref{conj always ULC} for the class of
negatively dependent measures.  

8.  Partial Symmetrization.  One could strengthen the preceding 
closure property by allowing symmetrization of only a subset
of the coordinates, for example, one could take $\mu' = (\mu +
\mu \circ \pi) / 2$ where $\pi$ is a transposition.  If one broadens
this to taking $\mu' = (1 - \ee) \mu + \ee \mu \circ \pi$, then
by iterating these with $\ee \rightarrow 0$, one obtains closure
under an arbitrary time-inhomogeneous stirring operation.  That 
is, let $\{ \pi_t : t \geq 0 \}$ be a $S_n$-valued stochastic Markov
process, with transitions from $\pi$ to $\tau \circ \pi$ at rates
$C(\tau , t)$ for each transposition $\tau$, where the functions
$C(\tau , t)$ are some arbitrary real functions.  Fix $T > 0$ and
let $\mu' = \mu \circ \pi_T$.  We require that our class of negatively
dependent measures, if it contains $\mu$, to contain any such $\mu'$.

One motivation for considering such a strong closure property is that
we expect it to hold when $\mu$ is a point mass, since then $\mu'$
is the state of an exclusion process at a fixed time.  It seems reasonable
that if the initial state is random, chosen from a negatively dependent 
measure $\mu$, then the state at time $T$ should still be negatively 
dependent.  Another plausibility argument is that going from $\mu$
to $(1 - \ee) \mu +  \ee \mu \circ \tau$ is akin to sampling without
replacement.  It is shown in Joag-Dev and Proschan (1983, example 3.2 (a))
that the values of samples drawn without replacement from a fixed 
(real-valued) population are negatively associated.  If the initial 
population is random with a negatively dependent law, this should 
still be true.  

9.  Truncation.  Given $\mu$ on $\B_n$, let $\mu'$ be $\mu$ 
conditioned on $a \leq \sum_j X_j \leq b$.  We say that $\mu'$
is the truncation of $\mu$ to $[a,b]$.  We may ask that our class
be closed under truncation.  This seems the least controversial when
$a = b$ and we are conditioning on the sum $\sum_j X_j$.  In fact,
Block, Savits and Shaked (1982) define a collection of random variables
$\{ X_1 , \ldots , X_n \}$ to satisfy Condition N if there
is some collection $\{ Y_1 , \ldots , Y_{n+1} \}$ of random variables
satisfying the positive lattice condition~(\ref{eq pos lattice}) and
some number $k$ such that the law of $\{ X_1 , \ldots , X_n \}$ 
is the law of $\{ Y_1 , \ldots \, Y_n \}$ conditioned on 
$\sum_{j=1}^{n+1} Y_j = k$.  They show that many examples of
negatively dependent measures from Karlin and Rinott (1980) can
be represented this way, and that this implies negative association.  
In fact, Joag-Dev and Proschan (1983, Theorem~2.6) show that if
any random variables $\{ X_e : e \in E \}$ with
law $\mu$ satisfy 
\begin{equation} \label{eq stoch2}
(\mu \| \sum_e X_e = k+1) \stoch (\mu \| \sum_e X_e = k) ,
\end{equation}
then $(\mu \| \sum_e X_e = a)$ is negatively associated; a
result of Efron (1965) is that~(\ref{eq stoch2}) holds when the 
real-valued variables $X_e$ have densities that are log concave,
which together with Joag-Dev and Proschan's result yields the Karlin
and Rinott result.  

Conditioning on an entire interval $[a,b]$ may seem less natural; it is 
a special case of the next closure operation.

10.  Rank rescaling.  Given a measure $\mu$ on $\B_n$ and a log-concave
sequence $q_0 , \ldots , q_n$, define the {\em rank rescaling of} $\mu$
by $\{ q_j \}$ to be the measure $\mu'$ given by 
$$\mu' (x) = {q_{|x|} \mu (x) \over \sum_{y \in \B_n} q_{|y|} \mu (y)} .$$
Here $|y|$ denotes the rank of $y$ in $\B_n$, that is, the number of
coordinates of $y$ that are 1.  When $q_j = 
\one_{[a,b]} (j)$, this reduces to truncation.  Another special case is 
$q_j = r^j$, which is the same as imposing a uniform external field. 
Rank rescaling may be too strong a closure property to demand, so we give two
plausibility arguments.  Firstly, observe that rank rescaling commutes
with external fields.  Thus when $\mu$ is a product Bernoulli
measure, the rank rescaling of $\mu$ by $\{ q_j \}$ is just an 
exchangeable ULC measure plus an external field, which we know to be
CNA+.  Secondly, Theorem~\ref{th closed class} below shows that
the closure of $\S_0$ under rank rescaling is still contained in 
the class JNRD+.  Unfortunately, since 
projections do not commute with rank rescaling, this class is not
closed under projections, so we do not know whether adding 
rank rescaling to the list of closure operations results in 
measures that are negatively associated.

A concrete application in which we would like to have these 
closure properties is the random forest.  Let $G$ be a graph
with $n$ vertices and edge set $E(G)$ and define the
{\em uniform random forest} $\eta : E(G) \to \{ 0 , 1 \}$ 
to be chosen uniformly among subsets of $E(G)$ with no cycles.  Thus we
generalize the well studied spanning tree model by allowing more than
one component.  Peter Winkler (personal communication) asks whether 
any negative dependence can be shown for this model.  Together with
closure under truncation, this would imply negative correlations in
{\em constrained random forests}, the simplest one of these being
when $\eta$ is chosen from acyclic edge sets with cardinality either 
$n-1$ or $n-2$.  There seems to be no negative correlation result
known even in this simple setting.

\subsection{Building a class of negatively dependent measures from the inside}
\label{ss closed class}

In this section we prove the following theorem, showing that asking for
closure under rank rescaling is reasonable.
\begin{th} \label{th closed class}
Let $\S$ be the smallest class of measures containing laws of single
Bernoulli random variables and closed under products, external fields
and rank rescaling.  Then every measure in $\S$ is JNRD+.  
\end{th}

The theorem is proved in several steps.  

\noindent{\ul{Step 1: Represent each $\mu$ in $\S$ by a tree.}}
Observe that external fields commute with products and rank rescaling.
Since an external field changes a Bernoulli variable into another Bernoulli,
all measures in $\S$ are built from Bernoulli laws by products and
rank rescaling.  Let $\T$ be a finite rooted tree, with each leaf $e$ labeled by
a Bernoulli law $\nu_e$, and each interior vertex $v$ labeled by 
a log-concave sequence $\{ q^{(v)}_j \}$, whose length is one more 
than the number of leaves below $v$.  Associate a measure $\mu^v$
to each interior vertex $v$ recursively, by letting $\mu^v$ be the
rank rescaling by $\{ q^{(v)}_j \}$ of the product of the
measures associated with the subtrees of $v$.  Then the above
observation implies that every measure in $\S$ is the measure
associated with the root of such a tree $\T$, so that if the
measure is the law of $\{ X(e) : e \in S \}$ then the set of
leaves of $\T$ is precisely $S$.  We may assume without 
loss of generality that every interior vertex of $\T$ has precisely 
two children.  We also note that
log-concavity is closed under convolution and pointwise products,
and thus by an easy induction the rank sequence for every measure
$\mu^v$ associated with any vertex $v$ of such a tree is log-concave.  

\noindent{\ul{Step 2: Use Lemma~\ref{lem Markov}.}}  For any vertex
$v$ of $\T$, define $Y_v$ to be the sum of $X_e$ over all leaves $e$
lying below $v$ (the root is at the top).  
Suppose $e$ and $f$ are two leaves of $\T$ and
let $v$ be their meeting vertex, that is, the lowest vertex of $\T$ having
both $e$ and $f$ as descendants.  Let $e = e_0 , e_1 , \ldots , e_k , v , f_l , 
\ldots , f_0 = f$ be the geodesic connecting $e$ and $f$ in $\T$.  I claim 
that the sequence $\{ Y_{e_0} , \ldots , Y_{e_k} , Y_{f_l} , \ldots , Y_{f_0} \}$
is Markov, and that each is stochastically increasing in the previous one, 
except that $Y_{f_l}$ is stochastically decreasing in $Y_{e_k}$.  
The conclusion of this step, which follows immediately 
 from Lemma~\ref{lem Markov} once the claims are established,
is that $X_e$ and $X_f$ are negatively correlated.  

Establishing the Markov property is a diagram chase.  Use the notation
$g \geq v$ to denote that the leaf $g$ is a descendant of the vertex $v$.
Slightly stronger than the Markov property is the fact that the
collection $\{ X_g : g \geq e_{j-1} \}$ and the collection $\{ X_g : 
g \ngeq e_j \}$ are independent given $Y_{e_j}$.  To see that this
independence property holds, write
$$ \mu (X_g = x_g : g \in S) = C \prod_{g \in E} \nu_g (x_g) 
   \prod_{v~~{\rm interior}} q^{(v)}_{y_v} \, ,$$
where $y_v := \sum_{g \geq v} x_g$.  Now observe that the only
factors in the product depending both on values $x_g$ for $g \geq e_{j-1}$
and for $g \ngeq e_j$ depend only on the total $y_{e_j}$, giving us
the desired conditional independence.  

\noindent{\ul{Step 3: Verify the part of the claim involving stochastic 
dependence.}}  We first record a simple lemma.
\begin{lem} \label{lem abc}
Let $\{ a_n \}$, $\{ b_n \}$, $\{ c_n \}$ be finite sequences of nonnegative
real numbers, with $a_i b_j c_{i+j}$ not identically zero.  Let $X$ and $Y$
be random variables such that
\begin{equation} \label{eq abc}
\P (X = i , Y = j) = K a_i b_j c_{i+j} 
\end{equation}
for some normalizing constant, $K$.  Then
\begin{quotation}
$(i)$ $X \uparrow (X+Y)$ if $b$ is log-concave; $Y \uparrow (X+Y)$ if $a$
is log-concave; \\

$(ii)$ $(X+Y) \uparrow X$ if $b$ is log-concave; $(X+Y) \uparrow Y$ if $a$
is log-concave; \\

$(iii)$ $X \downarrow Y$ if $c$ is log-concave; $Y \downarrow X$ if $c$
is log-concave; \\
\end{quotation}
\end{lem}

\noindent{\sc Proof:}  By symmetry it suffices to prove the first half of each 
statement.  We use the fact that if $\mu$ and $\nu$ are probability measures 
on the integers with $\mu (x) / \nu (x)$ increasing in $x$, then $\mu
\stoch \nu$.  

For statement~$(i)$, let $\mu$ be the conditional distribution
of $X$ given $X+Y = j$, and let $\nu$ be the conditional distribution of
$X$ given $X+Y = j+1$ (we deal only with the interval of values of
$j$ for which we are conditioning on events of positive probability).  
Then $\mu (x) = C a_x b_{j-x} c_j$ for some 
constant $C$, while $\nu (x) = C' a_x b_{j+1-x} c_{j+1}$ for some $C'$.
Hence $\mu (x) / \nu (x) = C'' b_{j-x} / b_{j+1-x}$, which is decreasing
in $x$ as long as $\{ b_j \}$ is log-concave.  Statements~$(ii)$ 
and~$(iii)$ are proved similarly.  For~$(iii)$, let $\mu$ be the 
conditional distribution of $X$ given $Y = j$ and $\nu$ be the 
conditional distribution of $X$ given $Y = j+1$.  Then
$\mu (x) / \nu (x) = C c_{j+x} / c_{j+x+1}$, which is increasing
in $x$ if $\{ c_j \}$ is log-concave.  And for~$(ii)$, let $\mu$
be the conditional distribution of $X + Y$ given $X = j$ and $\nu$
be the conditional distribution of $X + Y$ given $X = j+1$.  Then 
$\mu (x) / \nu (x) = C b_{x-j} / b_{x-j-1}$, which decreases
in $x$ when $\{ b_j \}$ is log-concave.   $\Cox$

The stochastic increases in the sequence $\{ Y_{e_0} , \ldots , Y_{e_k} , 
Y_{f_l} , \ldots , Y_{f_0} \}$ are now easy to verify.  Let $w$ be the
child of $e_{j+1}$ that is not $e_j$, let $X = Y_{e_j}$,
and let $Y = Y_w$.  Recall from the recursive construction
of the measures that $\mu^{e_j}$ gives $X$ a log-concave
sequence of probabilities, call it $\{ a_i \}$, that $\mu^w$ gives 
$Y$ a log-concave sequence of probabilities, call it $\{ b_i \}$, 
and that $\mu^{e_{j+1}}$ gives probabilities as in~(\ref{eq abc})
with $c_i = q^{e_{j+1}}_i$.  Replacing $\mu^{e_{j+1}}$ by the
measure $\mu$ associated with the root of the tree effectively 
alters the sequence $\{ c_i \}$ but not $\{ a_i \}$ or $\{ b_i \}$.
Since the sequences $\{ a_i \}$ and $\{ b_i \}$ are log-concave,
parts~$(i)$ and~$(ii)$ of the previous lemma imply that $X$ is
stochastically increasing in $X+Y$ and vice versa.  Since $X+Y = 
Y_{e_{j+1}}$, and since the argument works equally well for 
$f_j$ instead of $e_j$, this gives all parts of the claim except 
the fact that $Y_{f_l} \downarrow Y_{e_k}$.  

Let $v$ be the common parent of $e_k$ and $f_l$.  As before, 
we see that under the law $\mu^v$, $Y_{f_l}$ is stochastically
decreasing in $Y_{e_k}$, according to statement~$(iii)$ of
the lemma with $c_i = q^{(v)}_i$ which is log-concave.  Transferring
this argument to the measure $\mu$ is mostly a matter of using
the right notation to make it clear that the new sequence $\{ c_i \}$
is log-concave.  Let $v = v_0 , v_1 , \ldots , v_r$ be the path leading from $v$ 
to the root, and for $1 \leq i \leq r$, let $w_i$ be the child of $v_i$
not equal to $v_{i-1}$.  Let $a_i = \mu^{e_k} (Y_{e_k} = i)$ and
$b_i = \mu^{f_l} (Y_{f_l} = i)$.  Let $s^j_i = q^{(v_j)}_i$ and
let $t^j_i = \mu^{w_j} (Y_{w_j} = i)$.  Use the recursive definition of 
the measures $\mu^g$ to see that 
\begin{eqnarray*}
&& \mu (Y_{e_k} = i , Y_{f_l} = j) \\[1ex]
& = & K a_i b_j c_{i+j} \sum_{u_1, \ldots , u_r} \prod_{j=1}^r 
   t^j_{u_j} s^j_{i+j+u_1 + \cdots + u_j} \, .
\end{eqnarray*}
The summation term may be written as 
\begin{equation} \label{eq had}
((\cdots ((s^r * \overline{t^r}) \had s^{r-1}) * \overline{t^{r-1}} 
   \cdots \had s^1 ) * t^1) ,
\end{equation}
where $*$ denotes convolution, $\had$ denotes pointwise product, 
$\overline{~~~}$ denotes reversal, and $s^j$ and $t^j$ denote the sequences
$\{ s^j_i \} $ and $\{t^j_i \}$.  Since convolution, pointwise product and
reversal preserve log-concavity, this shows that the third part of
Lemma~\ref{lem abc} still applies, and finishes the verification.

\noindent{\ul{Step 4: Negative correlation implies h-NLC+.}}  
Observe that the property h-NLC+ is the same as NC+, where NC denotes
pairwise negative correlation.  To see this, note that an external
field with $W(e) \rightarrow 0$ or $\infty$ corresponds to conditioning
on $X_e = 0$ or 1 respectively.  Thus NC+ is equivalent to negative
correlation of any pair of variables, given values of any others, under
any external field, which is h-NLC+.  The conclusion of steps~2 and~3
were the NC property, and hence NC+, since the class is already closed
under external fields.  

\noindent{\ul{Step 5: Modifying the argument to get JNRD+.}}  
Let $e$ be a leaf of $\T$ and let $v_0 , v_1 , \ldots , v_k$ be the 
path from $e$ to the root, with $v_0 = e$.  Let $w_i$ be the child of
$v_i$ other than $v_{i-1}$.  I claim that the vector 
$(Y_{w_1} , \ldots , Y_{w_k})$ is stochastically 
decreasing in $X_e$.  This is shown by coupling, inducting
on $i$.  We will define a sequence $(Y_1 , \ldots , Y_k)$ to have
the conditional distribution of $(Y_{w_1} , \ldots , Y_{w_k})$
given $X_e = 0$ and $(Y_1' , \ldots , Y_k')$ to have
the conditional distribution of $(Y_{w_1} , \ldots , Y_{w_k})$
given $X_e = 1$ so that $(Y_1 - Y_1' , \ldots , Y_k - Y_k')$
has all coordinates zero except possibly for a single 1.  

When $i = 1$, we have $Y_{w_1} \downarrow X_e$
by part~$(iii)$ of Lemma~\ref{lem abc}, using log-concavity of
a sequence analogous to~(\ref{eq had}).  Since also $Y_{w_1} + 
X_e \uparrow X_e$ by part~$(i)$ of the lemma and log-concavity 
of the rank sequence for $Y_{w_1}$, this means we can define
$Y_1$ and $Y_1'$ so that $Y_1$ has the distribution
of $Y_{w_1}$ given $X_e = 0$, $Y_1'$ has the distribution
of $Y_{w_1}$ given $X_e = 1$ and $Y_1' + 1 \geq Y_1 \geq Y_1'$.
If $Y_1 = Y_1'+1$, then choose $(Y_2 , \ldots , Y_k)$ to have
the conditional distribution of $(Y_{w_2} , \ldots , Y_{w_k})$
given $X_e = 0$ and $Y_{w_1} = Y_1$.  This is the same as
the conditional distribution of $(Y_{w_2} , \ldots , Y_{w_k})$
given $X_e = 1$ and $Y_{w_1} = Y_1'$, so we may choose 
$(Y_2' , \ldots , Y_k') = (Y_2 , \ldots , Y_k)$.  If $Y_1 = Y_1'$,
then choose $Y_2$ and $Y_2'$ from the conditional distribution for
$Y_{w_2}$ given respectively that $Y_{v_1} = Y_1 +1$ and $Y_1$.
Again $Y_2' + 1 \geq Y_2 \geq Y_2'$, and we continue, setting
the remaining coordinates equal if $Y_2 = Y_2' + 1$, and otherwise
choosing $Y_3$ and $Y_3'$ and so on.  

The collections $\{ X_f : f \geq w_i \}$ are conditionally independent
as $i$ varies given $\{ Y_{w_i} : 1 \leq i \leq k \}$.  Thus we
may write the conditional law of $\{ X_f : f \neq e \}$ given $X_e = 0$ 
as a mixture over values $(r_1 , \ldots , r_k)$ of $(Y_1 , \ldots , Y_k)$
of product measures $\prod_{j=1}^k \mu_{j , r_j}$ where $\mu_{j,r_j}$
is the conditional law of $\{ X_e : e \geq w_j \}$ given $Y_{w_j} = r_j$.
The conditional law of $\{ X_f : f \neq e \}$ given $X_e = 1$ is the
same, but with a stochastically smaller mixing measure. 
Suppose the laws $\mu_{j,r_j}$ are stochastically increasing in $r_j$.
Then by stochastic comparison of the mixing measures, we see that
the conditional law of $\{ X_f : f \neq e \}$ given $X_e = 0$ 
dominates the conditional law of $\{ X_f : f \neq e \}$ given $X_e = 1$. 
The measures $\mu_{j,r_j}$ are in the class $\S$ ($\S$ is not closed
under projection but projections onto all variables in a subtree is OK).
Thus all that remains to verify JNRD+ is to prove the supposition, which
is the following lemma. 
\begin{lem} \label{lem Y stoch}
For any measure $\mu$ in the class $\S$, the conditional distribution
of $\mu$ given $\sum_e X_e = k+1$ stochastically dominates the
conditional distribution given $\sum_e X_e = k$.  
\end{lem}

To prove this we strengthen Lemma~\ref{lem abc} a little.  Recall that 
an element of a partially ordered set {\em covers} another if it is
greater and there is no element in between.  Say that a measure
$\mu$ on a partially ordered set covers the measure $\nu$ if there are 
random variables $X \sim \mu$ and $Y \sim \nu$ such that $X=Y$ or
$X$ covers $Y$.  
\begin{lem} \label{lem covers}
Under the hypotheses of Lemma~\ref{lem abc}, if $\{ a_n \}$ is log-concave, 
then $(X \| X+Y = k+1)$ covers $(X \| X+Y = k)$ and if $\{ c_n \}$ is
log-concave then $(X+Y \| X=k+1)$ covers $(X+Y \| X=k)$.
\end{lem}

\noindent{\sc Proof:} The likelihood ratio of the law of $X$ conditioned
on $X+Y = k+1$ to the law of $X+1$ conditioned on $X+Y=k$, evaluated
at the point $x$, is equal to $a_x b_{k+1-x} c_{k+1}/ (a_{x-1} b_{k+1-x}
c_k) = (c_{k+1}/c_k) (a_x / a_{x-1})$.  This is decreasing in $x$ by
log-concavity of $\{ a_n \}$.  The likelihood ratio of the law of 
$X+Y$ given $X=k+1$ to the law of $X+Y+1$ given $X=k$, evaluated
at the point $z$, is $a_{k+1} c_z / (a_k c_{z-1})$ which is decreasing
in $z$ by log-concavity of $\{ c_n \}$.   $\Cox$.

\noindent{\sc Proof of Lemma}~\ref{lem Y stoch}:  Induct on the height 
of the tree $\T$.  If $\T$
is a single leaf, then the statement is trivial.  Now suppose the root
of $\T$ has children $v$ and $w$ and assume for induction that
the lemma holds for $\mu^v$ and $\mu^w$.  Since the rank sequences
for $Y_v$ and $Y_w$ are log-concave, part~$(i)$
of Lemma~\ref{lem abc} show that $Y_v$ and $Y_w$ are each stochastically
increasing in $Y_v + Y_w$.  By Lemma~\ref{lem covers}, in fact the
law of $Y_v$ given $Y_v + Y_w = k+1$ covers the law of $Y_v$ given 
$Y_v + Y_w = k$, from which we conclude that the pair $(Y_v , Y_w)$ 
is stochastically increasing in $Y_v + Y_w$.  By the inductive hypothesis, 
$\{ X_e : e \geq v \}$ is stochastically increasing in $Y_v$ and the 
same is true with $v$ replaced by $w$.  Since $\{ X_e : e \geq v \}$ and 
$\{ X_e : e \geq w \}$ are conditionally independent given $Y_v$ and $Y_w$, 
this finishes the proof.   $\Cox$

\subsection{Further observations and conjectures}

Lemma~\ref{lem Y stoch} seems to be true in the following greater generality.
\begin{conj} \label{conj rank cover}
If $\mu$ is CNA+ then the conditional distribution $\mu$ given
$\sum_e X_e = k+1$ stochastically dominates the conditional
distribution $\mu$ given $\sum_e X_e = k$.  
\end{conj}
{\em Remark:} The conclusion of this conjecture appears in Joag-Dev
and Proschan (1983) as a hypothesis implying negative association.  
Does this condition fit into the theory of negative dependence better 
as a hypothesis or a conclusion?  The same could be asked about the ULC 
condition, c.f Conjectures~\ref{conj always ULC}~-~\ref{conj ULC implies NA}.

Another conjecture that seems to be true is as follows.  
\begin{conj} \label{conj ND cover}
If $\mu$ on $\B_n$ is CNA+ then the conditional distribution 
on $\B_{n-1}$ given $X_n = 0$ stochastically covers the conditional
distribution given $X_n = 1$.
\end{conj}

These conjectures may be strengthened by weakening the hypothesis
to JNRD+ or h-NLC+, but the + condition is essential, at least for
the second conjecture, as shown by the following example.

\noindent{\em Example:} Let $\mu$ be the measure on $\B_3$ with
equal probabilities $1/5$ for the points $(0,0,0), (0,0,1)$, $(0,1,0),
(1,0,0)$ and $(1,1,0)$.  This is CNA but not h-NLC+ (impose an external
field with $W(1)$ very small).  The measure $(\mu \| X_3 = 0)$ is
stochastically greater than the measure $(\mu \| X_3 = 1)$ but is
too much greater to cover it.

\begin{qq} 
Under what hypotheses on $\mu$ can one prove that 
\begin{equation} \label{eq cov}
(\mu \| \sum_e X_e = k+1) \stoch (\mu \| \sum_e X_e = k) ?
\end{equation}
\end{qq}
An answer to this question would be important for the following reason.
Let $A$ be any upset.  If we can establish~(\ref{eq cov}), then
$A \uparrow \sum_e X_e$ and in particular these have nonnegative 
covariance.  Therefore $A$ and $X_e$ have nonnegative covariance
for some $e$ and we have established proprty~$(ii)$ of Section~\ref{ss FM}. 
In particular, Conjecture~\ref{conj rank cover} implies
Conjecture~\ref{conj horizontal}.

\noindent{\bf Acknowledgements:}  Most of the blame for this goes
to Peter Doyle for egging me on in the early going and for proving 
Theorem~\ref{th closed class} with me.  Thanks to Peter Shor for 
suggestions pertaining to the urn model.  Thanks to Yosi Rinott for 
some helpful discussions on a previous draft of this paper.

\end{document}